\documentclass[12pt]{article}
\usepackage{amsfonts}
\usepackage{amsmath,amsthm,hyperref}
\usepackage[all]{xy}

\newtheorem{theorem}{Theorem}[section]
\newtheorem{proposition}[theorem]{Proposition}
\newtheorem{definition}[theorem]{Definition}
\newtheorem{corollary}[theorem]{Corollary}
\newtheorem{lemma}[theorem]{Lemma}

\theoremstyle{remark}
\newtheorem{remark}[theorem]{Remark}
\newtheorem{example}[theorem]{\bf Example}

\begin{document}
\title{\bf{Spacelike Willmore surfaces in 4-dimensional
Lorentzian space forms}}
\author{{Xiang Ma, Peng Wang
\footnote{ ~e-mail: $maxiang@math.pku.edu.cn,
wangpeng@math.pku.edu.cn$}}
\\{\small School of Mathematical Sciences,
 Peking University, 100871 }\\{\small People's Republic of China}}
\date{}
\maketitle

\begin{center}
{\bf Abstract}
\end{center}

Spacelike Willmore surfaces in 4-dimensional Lorentzian space
forms, a topic in Lorentzian conformal geometry which parallels
the theory of Willmore surfaces in $S^4$, are studied in this
paper. We define two kinds of transforms for such a surface, which
produce the so-called left/right polar surfaces and the adjoint
surfaces. These new surfaces are again conformal Willmore
surfaces. For them holds interesting duality theorem. As an
application spacelike Willmore 2-spheres are classified. Finally
we construct a family of homogeneous spacelike Willmore tori.\\

{\bf Keywords:} Spacelike Willmore surfaces; adjoint transforms;
polar surfaces; duality theorem

\section{Introduction}

Willmore surfaces are the critical surfaces with respect to the
conformally invariant Willmore functional. Many interesting results
related to them have been obtained (see [2,4,10,16]), and now they
are recognized as one of the most important surface classes in
M\"obius geometry.

For Lorentzian space forms there is also a parallel theory of
conformal geometry. Thus it is natural to generalize the notion of
Willmore surfaces to such a context. This idea was first followed by
Alias and Palmer in [1]. They considered the codim-1 case and
established such a theory as Bryant did in [2]: the conformal Gauss
map was introduced; the Willmore functional was defined as the area
with respect to the metric induced from this map; a surface is
Willmore if, and only if, its conformal Gauss map is harmonic. Later
Deng and Wang [8] treated timelike Willmore surfaces in Lorentzian
3-space; Nie [17] established a theory of conformal geometry about
hypersurfaces in Lorentzian space forms and computed the first
variation of Willmore functional.

In this paper we take the next step to study spacelike Willmore
surfaces in $Q^4_1$, the conformal compactification of the
4-dimensional Lorentzian space forms $R^{4}_{1},S^{4}_{1}$ and
$H^{4}_{1}$. In many aspects the theory is almost the same as in
M\"obius geometry, except that we have a distinctive construction as
below.

For a spacelike surface $[Y]$ immersed into $Q^4_1$, the normal
plane is Lorentzian at each point. The null lines $[L],[R]$ in this
plane define two conformal maps into $Q^4_1$, called \emph{the left
and the right polar surface}, while these transforms are called
($-$)transform and ($+$)transform, respectively. Conversely, $Y$ is
also the right polar surface of $[L]$, and the left polar surface of
$[R]$ (when $[L]$ and $[R]$ are immersions). That means
($-$)transform and ($+$)transform are mutual inverses to each other
(this is true even without the Willmore condition). Applying these
transforms successively, we obtain a sequence of conformal surfaces
as described by the following diagram:
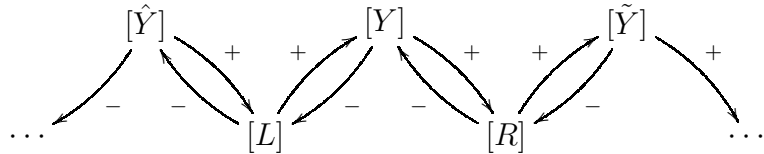
\[
\xymatrix{& &[\hat{Y}] \ar@/^/[dl]^{-} \ar@/^/[dr]^{+} & &[Y]
\ar@/^/[dl]^{-} \ar@/^/[dr]^{+}
& &[\tilde{Y}] \ar@/^/[dl]^{-} \ar@/^/[dr]^{+} &\\
&\cdots & &[L]\ar@/^/[ul]^{-} \ar@/^/[ur]^{+} & &[R]\ar@/^/[ul]^{-}
\ar@/^/[ur]^{+} & &\cdots}
\]
Our main result says that they are all Willmore surfaces if $[Y]$ is
assumed to be so.

It is interesting to notice that the two-step transforms
$[\hat{Y}]$ and $[\tilde{Y}]$ are located on \emph{the central
sphere} of the original Willmore surface $[Y]$ at corresponding
point, which mimics the property of the adjoint transforms in
$S^n$ as introduced by the first author [14] (indeed they could be
introduced in the same manner). In the special case that
$[\hat{Y}]=[\tilde{Y}]$, this yields a Willmore surface sharing
the same central sphere congruence as $[Y]$. It generalizes the
duality theorem of Bryant [2] and Ejiri [10], and such surfaces
will still be called \emph{S-Willmore} surfaces as in [10,13,14].
In particular, there is a surprising analogy between our
transforms and the so-called \emph{forward and backward B\"acklund
transforms} defined by Burstall et al. for Willmore surfaces in
$S^4$ [4].

When the underlying surface $M$ is compact, an important problem
is to classify all Willmore immersions of $M$ and to find the
values of their Willmore functionals (i.e. to determine the
critical values and critical points of the Willmore functional).
For Willmore 2-spheres in $S^3$ and $S^4$ this question was
perfectly answered by Bryant [2] and Montiel [14], respectively.
Precisely speaking, any Willmore 2-spheres in $S^4$ is the
conformal compactification of a complete minimal surface in
$R^{4}$, or the twistor projection of a complex curve in the
twistor space $\mathbb{C}P^{3}$. This follows from the duality
theorem and the vanishing theorem about holomorphic forms on
$S^2$. By the same method we could obtain similar characterization
result in the Lorentzian space.

\begin{flushleft} \textbf{Theorem}.
\begin{itshape}
Any spacelike Willmore 2-sphere in $Q^{4}_{1}$ is either the
conformal compactification of a complete spacelike stationary
surface (i.e. H=0) in $R^{4}_{1}$, or a polar surface of such a
surface (in the latter case the surface is the twistor projection
of a holomorphic curve in the twistor space of $Q^4_1$). For a
surface of the second type, its Willmore functional always equals
zero.
\end{itshape}
\end{flushleft}

This paper is organized as follows. In Section~2, we describe the
Lorentzian conformal space $Q^{4}_{1}$ as well as round 2-spheres in
it. The general theory about spacelike surfaces and the
characterization of Willmore surfaces are given in Section~3 and
Section~4. Then we study the transforms of spacelike Willmore
surfaces in Section~5. These transforms are utilized to classify
spacelike Willmore 2-spheres in Section~6. Finally we discuss some
special examples in Section~7 and construct a family of homogeneous
spacelike Willmore tori which are not S-Willmore.

In the sequel $y:M\to Q^4_1$ will always denote a smooth spacelike
immersion from an oriented surface $M$ unless it is explicitly
claimed otherwise.

\section{Lorentzian conformal geometry of $Q^{4}_{1}$}

Let $\mathbb{R}^n_s$ be the space $\mathbb{R}^n$ equipped with the
quadric form
\[
\langle x,x\rangle=\sum^{n-s}_{1}x^2_i-\sum^n_{n-s+1}x^{2}_i.
\]
In this paper we will mainly work with $\mathbb{R}^6_2$ whose light
cone is denoted as $C^5$. The quadric
\[
Q^4_1=\{\ [x]\in\mathbb{R}P^5\ |\ x\in C^5\setminus \{0\} \}
\]
is exactly the projectived light cone. The standard projection
$\pi:C^5\setminus\{0\}\rightarrow Q^4_1$ is a fiber bundle with
fiber $\mathbb{R}\setminus\{0\}$. It is easy to see that $Q^4_1$ is
equipped with a Lorentzian metric induced from projection $S^3\times
S^1\rightarrow Q^4_1$. Here
\begin{equation}S^3\times
S^1=\{x\in\mathbb{R}^6_2\ |\
\sum^4_{i=1}x^{2}_{i}=x^2_5+x^{2}_6=1\}\subset C^5\setminus\{0\}
\end{equation}
is endowed with the Lorentzian metric $g(S^3)\oplus (-g(S^1))$,
where $g(S^3)$ and $g(S^1)$ are standard metrics on $S^3$ and
$S^1$. So there is a conformal structure of Lorentzian metric
$[h]$ on $Q^4_1$. By a theorem of Cahen and Kerbrat [6], we know
that the conformal group of $(Q^4_1,[h])$ is exactly the
orthogonal group $O(4,2)/\{\pm1\}$, which keeps the inner product
of $\mathbb{R}^6_2$ invariant and acts on $Q^4_1$ by
\begin{equation}
T([x])=[xT],\ T\in O(4,2).
\end{equation}

As in the Riemannian case, there are three 4-dimensional
Lorentzian space forms, each with constant sectional curvature
$c=0,+1,-1$, respectively. They are defined by
\begin{align*}
R^4_1 &, ~c=0; \\
S^4_1 &:=\{x\in\mathbb{R}^5_1\ |\ \langle x,x\rangle=1 \}, c=1; \\
H^4_1 &:=\{x\in\mathbb{R}^5_2\ |\ \langle x,x\rangle=-1 \}, c=-1.
\end{align*}
Each of them could be embedded as a proper subset of $Q^4_1$:
\begin{equation}\label{varphi}
\begin{array}{llll}
\varphi_{0}:R^4_1\rightarrow Q^4_1, ~&
\varphi_{0}(x)=[(\frac{-1+\langle x,x\rangle}{2},
x,\frac{1+\langle x,x\rangle}{2})]; \\[1mm]
\varphi_{+}:S^4_1\rightarrow Q^4_1,
~& \varphi_{+}(x)=[(x,1)]; \\[1mm]
\varphi_{-}:H^4_1\rightarrow Q^4_1, ~& \varphi_{-}(x)=[(1,x)].
\end{array}
\end{equation}
It is easy to verify that these maps are conformal embeddings. In
particular, the Lorentzian conformal space $Q^4_1$ could be viewed
as the conformal compactification of $R^4_1$ by attaching the
light-cone at infinity to it, i.e.
\[Q^4_1=\varphi_{0}(R^4_1)\cup C_{\infty},\]
where $C_{\infty}=\{(a,u,a)\in\mathbb{R}P^5\ |\ \langle
u,u\rangle=0,a\in\mathbb{R}\}$. Thus $Q^4_1$ is the proper space to
study the conformal geometry of these Lorentzian space forms.

We note that the description above is valid in $n$-dimensional
space. The whole theory parallels M\"obius geometry, and Lorentzian
space forms are viewed as conic sections of $Q^n_1$.

Lorentzian conformal geometry is also analogous to M\"obius geometry
in that we have round spheres as the most important conformally
invariant objects. For our purpose here we only discuss \emph{round
2-spheres} (they were named \emph{conformal 2-spheres} in [1]). Each
of them could be identified with a 4-dim Lorentzian subspace in
$\mathbb{R}^6_2$. Given such a 4-space $V$, the round 2-sphere is
given by
\[S^2(V):=\{[v]\in Q^4_1 ~|~  v\in V\}.\]
Such spheres share the same properties as the round 2-spheres in
M\"obius geometry: they are not only topological 2-spheres, but also
geodesic 2-spheres when viewed as subsets of some Lorentzian space
form; they are totally umbilic spacelike surfaces. In our terms the
moduli space $\Sigma$ of all round 2-spheres in $Q^4_1$ can be
identified with the Grassmannian manifold
\[G_{3,1}(\mathbb{R}^{6}_{2}):=
\{\text{4-dim Lorentzian subspaces of}~\mathbb{R}^{6}_{2}\}.\]

\section{Basic equations for a surface in $Q^{4}_{1}$}

For a surface $y:M\rightarrow Q^{4}_{1}$ and any open subset
$U\subset M$, a local lift of $y$ is just a map $Y:U\rightarrow
C^5\setminus\{0\}$ such that $\pi\circ Y=y$. Two different local
lifts differ by a scaling, so the metric induced from them are
conformal to each other.

Let $M$ be a Riemann surface. An immersion $y:M\rightarrow
Q^{4}_{1}$ is called a conformal spacelike surface, if $\langle
Y_{z},Y_{z}\rangle=0$ and $\langle Y_{z},Y_{\bar{z}}\rangle >0$
for any local lift $Y$ and any complex coordinate $z$ on $M$.
(Here $Y_z=\frac{1}{2}(Y_u -i Y_v)$ is the complex tangent vector
for $z=u+iv$, and $Y_{\bar{z}}$ its complex conjugate.) For such a
surface there is a decomposition $M\times
\mathbb{R}^{6}_{2}=V\oplus V^{\perp}$, where
\begin{equation}
V={\rm Span}\{Y,{\rm d}Y,Y_{z\bar{z}}\}
\end{equation}
is a Lorentzian rank-4 subbundle independent to the choice of $Y$
and $z$. The orthogonal complement $V^{\perp}$ is also a
Lorentzian subbundle, which might be identified with the normal
bundle of $y$ in $Q^{4}_{1}$. Their complexifications are denoted
separately as $V_{\mathbb{C}}$ and $V^{\perp}_{\mathbb{C}}$.

Fix a local coordinate $z$. There is a local lift $Y$ satisfying
$|{\rm d}Y|^2=|{\rm d}z|^2$, called the canonical lift (with respect
to $z$). Choose a frame $\{Y,Y_{z},Y_{\bar{z}},N\}$ of
$V_{\mathbb{C}}$, where $N\in\Gamma(V)$ is uniquely determined by
\begin{equation}\label{eq-N}
\langle N,Y_{z}\rangle=\langle N,Y_{\bar{z}}\rangle=\langle
N,N\rangle=0,\langle N,Y\rangle=-1.
\end{equation}
For $V^{\perp}$ which is a Lorentzian plane at every point of $M$, a
natural frame is $\{L,R\}$ such that
\begin{equation}\label{eq-LR}
\langle L,L\rangle=\langle R,R\rangle=0,\langle L,R\rangle=-1.
\end{equation}
So $L$ and $R$ span the two null lines in $V^{\perp}$ separately.
They are determined up to a real factor around each point.

Given frames as above, it is straightforward to write down the
structure equations of $Y$. First note that $Y_{zz}$ is orthogonal
to $Y$, $Y_{z}$ and $Y_{\bar{z}}$. So there must be a complex
function $s$ and a section $\kappa\in
\Gamma(V_{\mathbb{C}}^{\perp})$ such that
\begin{equation}
Y_{zz}=-\frac{s}{2}Y+\kappa.
\end{equation}
This defines two basic invariants $\kappa$ and $s$ depending on
coordinates $z$. Similar to the case in M\"obius geometry, $\kappa$
and $s$ are interpreted as \emph{the conformal Hopf differential}
and \emph{the Schwarzian} of $y$, separately (see [5][14]).
Decompose $\kappa$ as
\begin{equation}\label{eq-kappa}
\kappa=\lambda_{1}L+\lambda_{2}R.
\end{equation}
Let $D$ denote the normal connection, i.e. the connection in the
bundle $V^{\perp}$. We have
\[D_z L=\alpha L,~~D_z R=-\alpha R \]
for the connection 1-form $\alpha{\rm d}z$. Denote
\begin{equation}\label{eq-Dkappa}
\langle \kappa,\bar\kappa\rangle=-\beta,~~ D_{\bar z}\kappa=\gamma_1
L +\gamma_2 R,
\end{equation}
where
\begin{equation}\label{eq-gamma}
\left\{\begin{array}{llll}
\beta=\lambda_1\bar{\lambda}_2+\lambda_2\bar{\lambda}_1,\\
\gamma_1=\lambda_{1\bar{z}}+\lambda_1\bar{\alpha},\\
\gamma_2=\lambda_{2\bar{z}}-\lambda_2\bar{\alpha}.
\end{array}\right.
\end{equation}
The structure equations are given as follows:
\begin{equation}\label{eq-moving}
\left\{\begin {array}{lllll}
Y_{zz}=-\frac{s}{2}Y+\lambda_{1}L+\lambda_{2}R,\\
Y_{z\bar{z}}=\beta Y+\frac{1}{2}N,\\
N_{z}=2\beta Y_{z}-sY_{\bar{z}}+2\gamma_{1}L+2\gamma_{2}R,\\
L_{z}=\alpha L-2\gamma_{2}Y+2\lambda_{2}Y_{\bar{z}},\\
R_{z}=-\alpha R-2\gamma_{1}Y+2\lambda_{1}Y_{\bar{z}},
\end {array}\right.
\end{equation}
The conformal Gauss, Codazzi and Ricci equations as integrable
conditions are:
\begin{equation}\label{eq-integ}
\left\{\begin {array}{lllll}
s_{\bar{z}}=-2\beta_{z}-4\lambda_{1}\bar{\gamma}_{2}
-4\lambda_{2}\bar{\gamma}_{1},\\
{\rm Im}(\gamma_{1\bar{z}}+\gamma_{1}\bar{\alpha}
+\frac{\bar{s}}{2}\lambda_{1})=0,\\
{\rm Im}(\gamma_{2\bar{z}}-\gamma_{2}\bar{\alpha}
+\frac{\bar{s}}{2}\lambda_{2})=0,\\
D_{\bar{z}}D_{z}L-D_{z}D_{\bar{z}}L
=2(\lambda_{2}\bar{\lambda}_{1}-\bar{\lambda}_{2}\lambda_{1})L,\\
D_{\bar{z}}D_{z}R-D_{z}D_{\bar{z}}R
=-2(\lambda_{2}\bar{\lambda}_{1}-\bar{\lambda}_{2}\lambda_{1})R.
\end {array}\right.
\end{equation}
These are quite similar to the theory in [5]. In particular, the
second and the third equation above could be combined and written as
a single conformal Codazzi equation:
\begin{equation}\label{eq-codazzi}
{\rm Im}(D_{\bar{z}}D_{\bar{z}}\kappa+\frac{\bar{s}}{2}\kappa)=0.
\end{equation}
\begin{remark}\label{rem:kappa}
Another important fact we will need later is that $\kappa~({\rm
d}z)^\frac{3}{2}({\rm d}\bar{z})^{-\frac12}$ is a globally defined
vector-valued complex differential form.
\end{remark}

\section{Willmore functional and Willmore surfaces}

\begin{definition}\label{def:1}
For a conformal spacelike surface $y:M\rightarrow Q^{4}_{1}$, the
4-dim Lorentzian subspace
\[V={\rm Span}\{Y,{\rm d}Y,Y_{z\bar{z}}\}\]
at one point $p\in M$ is identified with a round 2-sphere $S^2(V)$
in $Q^4_1$ as in Section~2. We call it \emph{the central sphere} of
the surface $y$ at $p$.
\end{definition}

The notion of central spheres comes from M\"obius geometry, where it
is of great importance in the study of surfaces (and general
submanifolds) [2,5,8,19]. It is also known as \emph{the mean
curvature sphere} of the immersed surface $y$ at $p$, characterized
as the unique round 2-sphere $y^{*}$ tangent to $y$ at $p$ and
sharing the same mean curvature vector as $y$ at this point. (The
ambient space is endowed with a metric of some space form). In the
Lorentzian case this is also true.

\begin{proposition}
A surface immersed in a lorentzian space form envelops its central
sphere congruence and shares the same mean curvature with these
round 2-spheres at corresponding points.
\end{proposition}
\begin{proof}
To prove the conclusion for the flat space $R^4_1$, consider a
surface $y:M\to R^4_1$ and a point $p\in M$. Let $y^{*}:S^2\to
R^4_1$ be the \emph{mean curvature sphere} associated with $y$ at
$p$ as characterized above. It suffices to show that $y^{*}$
coincides with the central sphere of $y$ at $p$. Embed surface $y$
into $Q^4_1$ via $\varphi_{0}$ as given by \eqref{varphi}, with lift
\[
Y =\Big(\frac{-1+\langle y,y\rangle}{2}, y,\frac{1+\langle
y,y\rangle}{2}\Big).
\]
Computation shows that the central sphere of $y$ at $p$, identified
with $V={\rm Span}\{Y,{\rm Re}(Y_z),{\rm Im}(Y_z),Y_{z\bar{z}}\}$,
is determined by the position vector $y(p)$, the tangent plane ${\rm
d}y(TM_p)$, and the mean curvature vector $H(p)$, which coincide
with those of $y^*$ by our assumption. So $y,y^*$ share the same
central sphere at $p$. Yet for the round 2-sphere $y^*$, its central
sphere at any point is exactly itself (they fall into the same 4-dim
subspace), which verifies our assertion. For surfaces in $S^4_1$ or
$H^4_1$ the proof is similar.
\end{proof}

\begin{corollary}\label{cor:mini}
In particular, if the central sphere congruence of $y$ is a family
of planes in $R^4_1$, this surface must have mean curvature zero at
every point, thus be a stationary surface in $R^4_1$.
\end{corollary}

The central sphere congruence is conformally invariant in the
sense that for two surfaces $[Y^{'}],[Y]$ differing to each other
by the action of $T\in O(4,2)$, their central spheres at
corresponding points also differ by this transformation. This
tells us that although the mean curvature sphere of a surface at
one point is defined in terms of metric geometry, it is indeed a
conformal invariant by Proposition~4.2 and the observation above.
Viewed as a map from $M$ to $\Sigma$, the moduli space of round
2-spheres, it has another name, \emph{the conformal Gauss map} of
$y$ [1,2]. In Section~2 we have identified $\Sigma$ with the
Grassmannian $G_{3,1}(\mathbb{R}^{6}_{2})$, which could be further
embedded into the space of multi-vectors (of certain type and of
length 1) in $\mathbb{R}^{6}_{2}$:
\[
\Sigma \simeq G_{3,1}(\mathbb{R}^{6}_{2})\hookrightarrow
\Lambda_{3,1}(\mathbb{R}^{6}_{2}).
\]
The latter is endowed with the canonical semi-Riemannian metric as
usual. This provides the appropriate framework for the discussion of
the geometry of the conformal Gauss map.

\begin{definition}
For a conformally immersed surface $y:M\to Q^4_1$ with canonical
lift $Y$ (with respect to a local coordinate $z$), define
\[
G:=Y\wedge Y_{u}\wedge Y_{v}\wedge N=-2i\cdot Y\wedge Y_{z}\wedge
Y_{\bar{z}} \wedge N,\ z=u+iv,
\]
where $N\equiv 2Y_{z\bar{z}}(\!\!\mod Y)$ is the frame vector
determined in \eqref{eq-N}. Note that $\langle G,G\rangle=1$ and
that $G$ is well defined. We call $G:M\rightarrow
G_{3,1}(\mathbb{R}^{6}_{2})$ \emph{the conformal Gauss map} of $y$.
It is noteworthy that $V\in G_{3,1}(\mathbb{R}^{6}_{2})$ determines
$V^{\perp}\in G_{1,1}(\mathbb{R}^{6}_{2})$ and vice versa. Hence the
geometry of $G$ is equivalent to the geometry of the associated map
\[
G^{\perp}:=L\wedge R: M\rightarrow G_{1,1}(\mathbb{R}^{6}_{2}),
\]
where $L,R$ are normal vectors as given in \eqref{eq-LR}.
\end{definition}

The conformal Gauss map is important in that it induces a
conformally invariant conformal metric. Direct computation using
\eqref{eq-moving} shows
\begin{proposition}
For a conformal surface $y:M\rightarrow Q^{4}_{1}$, $G$ induces a
metric
\[
g:=\frac{1}{4}\langle {\rm d}G,{\rm d}G\rangle=\langle
\kappa,\bar{\kappa}\rangle|dz|^{2}
\]
on M, where $\kappa=\lambda_1 L+\lambda_2 R$ is the conformal Hopf
differential. This metric might be positive definite, negative
definite, or degenerate according to the sign of
$\langle\kappa,\bar{\kappa}\rangle=-\beta=-(
 \lambda_{1}\bar{\lambda}_{2}+\lambda_{2}\bar{\lambda}_{1}).$
\end{proposition}
Now we can introduce the Willmore functional and Willmore surfaces.
\begin{definition}
\emph{The Willmore functional} of $y$ is defined as the area of M
with respect to the metric above:
\begin{equation}
W(y):=\frac{i}{2}\int_{M}|\kappa|^{2}dz\wedge d\bar{z}.
\end{equation}
An immersed surface $y:M\rightarrow Q^{4}_{1}$ is called a
\emph{Willmore surface}, if it is a critical surface of the
Willmore functional with respect to any variation of the map
$y:M\rightarrow Q^{4}_{1}$.
\end{definition}

Willmore surfaces can be characterized as follows, which is similar
to the conclusions in codim-1 case [1,8] as well as in M\"obius
geometry [2,5,10,14].

\begin{theorem}\label{thm-willmore}
For a conformal spacelike surface $y:M\rightarrow Q^{4}_{1}$, the
following four conditions are equivalent:

(i) $y$ is Willmore.

(ii) The conformal Gauss map $G$ is a harmonic map into
$G_{3,1}(\mathbb{R}^{6}_{2})$.

(iii) The conformal Hopf differential $\kappa$ of $y$ satisfies the
Willmore condition as below, which is stronger than the conformal
Codazzi equation \eqref{eq-codazzi}:
\begin{equation}\label{eq-willmore}
D_{\bar{z}}D_{\bar{z}}\kappa+\frac{\bar{s}}{2}\kappa=0.
\end{equation}

(iv) In a Lorentzian space form of sectional curvature $c$, $y$
satisfies the Euler-Lagrange equation
\begin{equation}\label{eq-willmore2}
\Delta \vec{H}-2(|\vec{H}|^2+K-c)\vec{H}=0.
\end{equation}
Here $\Delta, \vec{H}, K$ are the Laplacian of the induced metric,
the mean curvature vector, and the Gaussian curvature of $y$,
respectively.
\end{theorem}

The proof to Theorem~\ref{thm-willmore} is completely the same as in
M\"obius geometry (we refer the reader to [15]). Note that when we
take a variation $y_t$ of the immersion $y$, generally $y_t$ is not
conformal to $y$, hence we have to consider the variation of the
Willmore functional with respect to a varied complex structure $J_t$
over $M$. Yet one can verify that this change of complex structure
$J$ contributes nothing to the first variation of the Willmore
functional. Then the Willmore condition \eqref{eq-willmore} can be
derived easily.

The equivalent condition \eqref{eq-willmore2} in this theorem also
implies that stationary surfaces (i.e. surfaces with mean curvature
$\overrightarrow{H}=0$) in Lorentzian space forms are Willmore.
Indeed they belong to a subclass of Willmore surfaces, the so-called
\emph{S-Willmore surfaces}. The latter are exactly those Willmore
surfaces with dual surfaces (see the next section). See Ejiri [10]
and Ma [14] for the counterpart in M\"obius geometry.

\begin{definition}
A conformal Willmore surface $y:M\rightarrow Q^{4}_{1}$ is called
a \emph{S-Willmore surface} if it satisfies
$D_{\bar{z}}\kappa\,\|\, \kappa,$ i.e.,
$D_{\bar{z}}\kappa=-\frac{\bar{\mu}}{2}\kappa$ for some local
function $\mu$ when $\kappa\ne 0$.
\end{definition}

\begin{definition}
Let $y:M\rightarrow Q^{4}_{1}$ be a spacelike surface. We call $y$
\emph{null-umbilic} if its Hopf differential is isotropic, i.e.
$\langle\kappa,\kappa\rangle=0$ (equivalently, $\lambda_1$ or
$\lambda_2$ vanishes). $y$ is umbilic if $\kappa=0$ (equivalently,
$\lambda_1=\lambda_2=0$).
\end{definition}

So far our notions, constructions and results can all be easily
generalized to n-dimensional spaces. Yet in $Q^{4}_{1}$ null-umbilic
surfaces have a special meaning. They are equivalent to $\mathcal
{O}^{+}$-holomorphic maps into the twistor space of $Q^{4}_{1}$
according to [12]. Thus they are similar to isotropic surfaces in
$S^{4}$ (which are also twistor projection of complex curves). Yet
there are also important differences. For example, isotropic
surfaces in $S^{4}$ are always S-Willmore, yet for null-umbilic
surfaces this is not necessarily true. (Only under the additional
Willmore condition can we show that a null-umbilic surface is
S-Willmore.)

\section{Transforms of Spacelike Willmore surfaces}
In this section, we will define two transforms for surfaces in
$Q^{4}_{1}$ and show that the new surfaces derived from them are
also Willmore if the original surface is Willmore.

\subsection{Right/left polar surfaces; ($+$/$-$)transforms}

For a conformal spacelike surface $y:M\rightarrow Q^{4}_{1}$ with
canonical lift $Y:M\rightarrow R^6_2$ with respect to complex
coordinate $z=u+iv$, its normal plane at any point is spanned by
two lightlike vectors $L,R$. Suppose that $R^6_2$ is endowed with
a fixed orientation and that
\[
\{Y,Y_u,Y_v,N,R,L\}
\]
form a positively oriented frame. $\{R,L\}$ might also be viewed
as a frame of the normal plane compatible with the orientation of
$M$ and that of the ambient space. Since $\langle L,R\rangle=-1$
has been fixed in \eqref{eq-LR}, either one of the null lines
$[L]$ ($[R]$) is well-defined.
\begin{definition}
The two maps
\[ [L],[R]:M \rightarrow Q^{4}_{1} \]
are named {\bf the left and the right polar surface} of $y=[Y]$,
respectively.
\end{definition}
\begin{remark}
Denote
$e_{+}=\frac{1}{\sqrt{2}}(R-L),~e_{-}=\frac{1}{\sqrt{2}}(R+L).$ Then
$\{e_{+},e_{-}\}$ is a positively oriented orthonormal frame of the
normal plane, and $L,R$ could be written as
\[
L=\frac{1}{\sqrt{2}}(e_{-}-e_{+}),
~~R=\frac{1}{\sqrt{2}}(e_{-}+e_{+}).
\]
Thus we also call $[L]$ the {\bf ($-$)transform}, and $[R]$ the {\bf
($+$)transform} of $[Y]$. At the same time these names correspond to
the directions of these transforms in the diagram below:
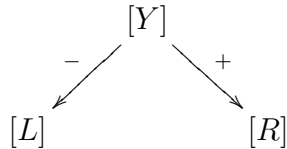
\[
\xymatrix{ &[Y] \ar[dl]_{-} \ar[dr]^{+} & \\
[L] & & [R]}
\]
The name \emph{polar surfaces} comes from Lawson's similar
construction for minimal surfaces in $S^3$ [11].
\end{remark}
\begin{proposition}\label{prop-inverse}
The polar surfaces $[L],[R]:M\to Q^{4}_{1}$ are both conformal
maps. $[L]$ ($[R]$) is degenerate if, and only if, $\lambda_2=0$
($\lambda_1=0$); it is a spacelike immersion otherwise. The
original surface $[Y]$ is the left polar surface of $[R]$ (the
right polar surface of $[L]$) when $[R]$ ($[L]$) is not
degenerate.
\end{proposition}
\begin{proof}
The first two conclusions for $[L]$ follow directly from
\[L_{z}=\alpha L-2\gamma_{2}Y+2\lambda_{2}Y_{\bar{z}}\]
by \eqref{eq-moving}. Differentiating this equation once more, by
\eqref{eq-gamma}\eqref{eq-moving} we find
\begin{align}
L_{z\bar{z}}=\bar{\alpha}L_z +&\alpha L_{\bar{z}}
+(\alpha_{\bar{z}}-\alpha\bar{\alpha}+2\lambda_2\bar{\lambda}_1)L
\notag \\
&+2\lambda_2\bar{\lambda}_2 \cdot R
-2\left(\gamma_{2\bar{z}}-\gamma_{2}\bar{\alpha}
+\frac{\bar{s}}{2}\lambda_{2}\right)Y.
\end{align}
When $\lambda_2\ne 0$, we can verify directly that $Y$ and
\begin{equation}
\hat{Y}=2\left|\frac{\gamma_2}{\lambda_2}\right|^2 Y
-2\frac{\bar\gamma_2}{\bar\lambda_2}~ Y_z
-2\frac{\gamma_2}{\lambda_2}~ Y_{\bar{z}} +N
+\frac{\gamma_{2\bar{z}}-\gamma_{2}\bar{\alpha}
+\frac{\bar{s}}{2}\lambda_{2}}{\lambda_2\bar\lambda_2}~ L
\end{equation}
are two lightlike vectors in the orthogonal complement of ${\rm
Span}\{L,L_z,L_{\bar{z}},L_{z\bar{z}}\}$ with $\langle
Y,\hat{Y}\rangle=-1$, and that
$\{L,L_z,L_{\bar{z}},L_{z\bar{z}},Y,\hat{Y}\}$ is again a positively
oriented frame. So for $[L]$, the left polar surface of $[Y]$, its
right polar surface is exactly $[Y]$. For $[R]$ the proof is
similar. In other words, the $(+)$transform is the inverse to the
$(-)$transform and vice versa when all surfaces concerned are
immersed.
\end{proof}
\begin{remark}
On the other hand, $[\hat{Y}]$ might be viewed as the 2-step
$(-)$transform of $y=[Y]$. Similarly we have the 2-step
$(+)$transform $[\tilde{Y}]$ as the right polar surface of $[R]$:
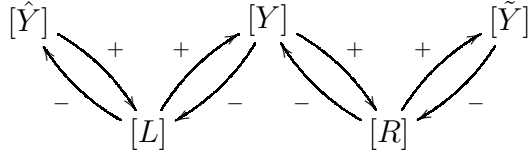
\[
\xymatrix{ &[\hat{Y}]  \ar@/^/[dr]^{+} & &[Y] \ar@/^/[dl]^{-}
\ar@/^/[dr]^{+}
& &[\tilde{Y}] \ar@/^/[dl]^{-} \\
& &[L]\ar@/^/[ul]^{-} \ar@/^/[ur]^{+} & &[R]\ar@/^/[ul]^{-}
\ar@/^/[ur]^{+} & }
\]
Note that $[L],[R]$ are also 2-step transforms to each other.
\end{remark}

\subsection{($+$/$-$)transforms preserve Willmore property}

Assume $y:M\rightarrow Q^{4}_{1}$ is an immersed spacelike
Willmore surface with canonical lift $Y:M\rightarrow R^6_2$ for a
given coordinate $z$ locally. We want to show that the
$(+)$transform and $(-)$transform again produce Willmore surfaces.

Assume that the left polar surface $[L]$ is an immersion, i.e.
$\lambda_{2}\ne 0$. Set
\begin{equation}\label{eq-mu}
-\frac{\bar\mu}{2}:=\frac{\gamma_2}{\lambda_2}.
\end{equation}
According to the conclusions of Theorem~\ref{thm-willmore} and
Proposition~\ref{prop-inverse}, we need to show that the conformal
Gauss map of $[L]$, represented by $Y\wedge\hat{Y}$, is a harmonic
map. The Willmore condition \eqref{eq-willmore} for $y$ amounts to
say
\begin{equation}\label{eq-willmore3}
\left\{\begin{array}{ll} \gamma_{1\bar{z}}+\gamma_{1}\bar{\alpha}
+\frac{\bar{s}}{2}\lambda_{1}=0,\\
\gamma_{2\bar{z}}-\gamma_{2}\bar{\alpha}
+\frac{\bar{s}}{2}\lambda_{2}=0.
\end{array}\right.
\end{equation}
Hence the expression of $\hat{Y}$ in (18) is simplified to
\begin{equation}\label{eq-yhat}
\hat{Y}:=\frac{|\mu|^2}{2} Y+\bar{\mu}Y_{z}+\mu Y_{\bar{z}}+N.
\end{equation}
The Willmore condition also implies
\begin{equation}\label{eq-muz}
\mu_{z}-\frac{\mu^{2}}{2}-s=0,
\end{equation}
because one can verify directly that
\[
2(\gamma_{2\bar{z}}-\gamma_{2}\bar{\alpha}
+\frac{\lambda_{2}\bar{s}}{2})
=(-\bar{\mu}_{\bar{z}}+\frac{\bar{\mu}^{2}}{2}+\bar{s})\lambda_{2}
\]
using the expressions of $\gamma_2$ \eqref{eq-gamma} and $\mu$
\eqref{eq-mu}.

For convenience of computation, set a new frame
\[ \{Y,\hat{Y},P,\bar{P},L,R\},
~~\text{with}~P:=Y_z+\frac{\mu}{2}Y, \] so that $\langle
Y,P\rangle=\langle \hat{Y},P\rangle=0.$ Differentiating $\hat{Y}$
and invoking \eqref{eq-muz}, we find
\begin{equation}\label{eq-Dyhat}
\hat{Y}_{z}=\frac{\mu}{2}\hat{Y}+ \rho P +\sigma L,
\end{equation}
where
\begin{equation}\label{eq-rho}
\left\{\begin {array}{ll}
\rho=\bar\mu_z+2\lambda_1\bar{\lambda}_2+2\lambda_2\bar{\lambda}_1
=\bar{\mu}_{z}+2\beta,\\
\sigma=2\gamma_{1}+\lambda_{1}\bar{\mu}.
\end {array}\right.
\end{equation}
For the frame $\{Y,\hat{Y},P,\bar{P},L,R\}$, the structure equations
are
\begin{equation}\label{eq-moving2}
\left\{\begin {array}{llllll}
Y_{z}=-\frac{\mu}{2}Y+P,\\[1mm]
\hat{Y}_{z}=~~\frac{\mu}{2}\hat{Y}+\rho P+\sigma L,\\[1mm]
P_{z}=~~\frac{\mu}{2}P+\lambda_{1}L+\lambda_{2}R,\\[1mm]
\bar{P}_{z}=-\frac{\mu}{2}\bar{P}
+\frac{\rho}{2}Y+\frac{1}{2}\hat{Y},\\[1mm]
L_{z}=~~\alpha L+2\lambda_{2}\bar{P},\\[1mm]
R_{z}=-\alpha R+2\lambda_{1}\bar{P}-\sigma Y.
\end {array}\right.
\end{equation}
Now the Willmore condition, \eqref{eq-muz} and the first one in
\eqref{eq-willmore3}, yields
\begin{align}
\rho_{\bar{z}} &=\bar{\mu}\rho-2\bar{\lambda}_{2}\sigma,
\label{eq-Drho}\\
\sigma_{\bar{z}}&=(-\bar{\alpha}+\frac{\bar{\mu}}{2})\sigma.
\end{align}
The computation is straightforward by the expressions
\eqref{eq-rho}\eqref{eq-gamma} and the first equation in
\eqref{eq-integ} (the conformal Gauss equation).

After these preparations, now we can compute out that
\begin{align*}
(Y\wedge\hat{Y})_{z}~=&~P\wedge\hat{Y}
+\rho Y\wedge P+\sigma Y\wedge L, \\
(Y\wedge\hat{Y})_{z\bar{z}} =&~\frac{\rho+\bar{\rho}}{2}
Y\wedge\hat{Y}+\sigma\bar{P}\wedge L+\bar{\sigma}P\wedge
L+\rho\bar{P}\wedge P+\bar{\rho}P\wedge\bar{P}\\
&~~+\underbrace{(\rho_{\bar{z}}-\bar{\mu}\rho
+2\sigma\bar{\lambda}_{2})}_{=0}Y\wedge
P+\underbrace{(\sigma_{\bar{z}}-\tfrac{\bar{\mu}}{2}\sigma
+\bar{\alpha}\sigma)}_{=0}Y\wedge L .
\end{align*}
Thus $Y\wedge\hat{Y}$ is a (conformal) harmonic map into
$G_{1,1}(\mathbb{R}^6_2)$ as desired. This shows that $[L]$ is
Willmore. For $[R]$ the proof is similar. Sum together, we have
proved
\begin{theorem}
Let $y:M\rightarrow Q^{4}_{1}$ be a spacelike Willmore surface. Then
its left and right polar surfaces $[L],[R]:M\rightarrow Q^{4}_{1}$
are also spacelike Willmore surfaces when they are not degenerate.
\end{theorem}
Since the $(-)$transform and the $(+)$transform preserve the
Willmore property, the same holds true for the 2-step transforms
$[\hat{Y}],[\tilde{Y}]$ when they are defined.
\begin{definition}
When $[Y]$ is Willmore, $[\hat{Y}],[\tilde{Y}]$ are also conformal
spacelike Willmore surfaces, called separately \emph{the left
adjoint transform} and \emph{the right adjoint transform} of $[Y]$.
\end{definition}
\begin{remark}\label{rem:adjoint}
Another equivalent way to define adjoint transforms of a given
Willmore surface is to follow the idea in [14]. In particular, the
adjoint transforms defined at here share many properties as before.
Taking $[\hat{Y}]$ for example, we have:
\begin{itemize}
\item[(1)] The (left) adjoint transform $[\hat{Y}]$ is conformal
to $[Y]$; it locates on the central sphere congruence of $[Y]$
according to \eqref{eq-yhat} and Definition~\ref{def:1}.
\item[(2)] $Y\wedge\hat{Y}$ is a conformal harmonic map into the
Grassmannian $G_{1,1}(\mathbb{R}^6_2)$.
\item[(3)] When the two adjoint transforms coincide, this surface
$[\hat{Y}]=[\tilde{Y}]$ will share the same central sphere
congruence with $[Y]$. (See the duality theorem in the next
subsection.)
\end{itemize}
The interested reader may confer [14] for a comparison. Here we
derive them from the polar surfaces, which seems more natural in our
context. Note that $[L],[R]$ are also adjoint transforms to each
other, as visualized below:
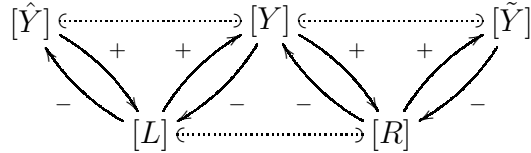
\[
\xymatrix{ &[\hat{Y}] \ar@/^/[dr]^{+} \ar@{.)}[rr]& &[Y]
\ar@/^/[dl]^{-} \ar@/^/[dr]^{+} \ar@{.)}[ll]\ar@{.)}[rr]
& &[\tilde{Y}] \ar@/^/[dl]^{-} \ar@{.)}[ll]\\
& &[L]\ar@/^/[ul]^{-} \ar@/^/[ur]^{+} \ar@{.)}[rr]&
&[R]\ar@/^/[ul]^{-} \ar@/^/[ur]^{+} \ar@{.)}[ll]& }
\]
The chain of $(-)$transforms and $(+)$transforms also demonstrates a
striking similarity with the backward and forward B\"acklund
transforms introduced for Willmore surfaces in $S^4$ [4]. In
particular, the 2-step B\"acklund transforms there could also be
identified with the adjoint transforms in [14]. An interesting
difference is that our $(-/+)$transforms are defined in a
conformally invariant way, whereas the 1-step B\"acklund transforms
are only properly defined in some affine space $R^4$.
\end{remark}

\subsection{Duality theorem of S-Willmore surfaces}

In the picture given above, a special case is noteworthy, namely
that when $[\hat{Y}]=[\tilde{Y}]$. This might be characterized by
the following
\begin{theorem}[Duality Theorem]
Let $y=[Y]:M\rightarrow Q^{4}_{1}$ be a spacelike S-Willmore surface
with polar surfaces $[L],[R]$ and adjoint transforms
$[\hat{Y}],[\tilde{Y}]$. Suppose that both of $[L],[R]$ are not
degenerate, i.e., $\lambda_1\ne 0, \lambda_2\ne 0$. Then the
conditions below are equivalent:
\begin{itemize}
\item[(1)] $[\hat{Y}]=[\tilde{Y}]$, i.e., the two adjoint
transforms coincide.
\item[(2)] $y=[Y]$ is a S-Willmore surface, i.e.
$D_{\bar{z}}\kappa=-\frac{\bar{\mu}}{2}\kappa$ for some $\mu$.
\item[(3)] $[\hat{Y}]$ (or $[\tilde{Y}]$) shares the same central
sphere congruence with $[Y]$.
\end{itemize}
\end{theorem}
\begin{proof}
When $y$ is Willmore, its right adjoint transform $[\tilde{Y}]$
might be given in a formula similar to \eqref{eq-yhat} with
\[
\tilde{Y}:=\frac{|\mu_1|^2}{2} Y+\bar{\mu}_1 Y_{z}+\mu_1
Y_{\bar{z}}+N,~~-\frac{\bar{\mu}_1}{2}:=\frac{\gamma_1}{\lambda_1}.
\]
Thus it is obvious that $[\hat{Y}]=[\tilde{Y}]$ if and only if
$-\bar{\mu}/2=\gamma_1/\lambda_1=\gamma_2/\lambda_2$, which is
equivalent to the S-Willmore condition. This shows
$"(1)\Leftrightarrow(2)"$. By \eqref{eq-moving2} we also know that
\[
{\rm Span}\{\hat{Y},\hat{Y}_z,\hat{Y}_{\bar{z}},\hat{Y}_{z\bar{z}}\}
={\rm Span}\{Y,Y_z,Y_{\bar{z}},Y_{z\bar{z}}\}
\]
if and only if $\sigma:=2\gamma_{1}+\lambda_{1}\bar{\mu}=0$, where
$-\bar{\mu}/2:=\gamma_2/\lambda_2$. This shows that $[\hat{Y}]$ has
the same central sphere congruence as $[Y]$ exactly when
$-\bar{\mu}/2=\gamma_1/\lambda_1=\gamma_2/\lambda_2$. So
$"(3)\Leftrightarrow(2)"$, and the proof is completed.
\end{proof}
\begin{remark}
Condition (3) in this theorem tells us that when $[Y]$ is
S-Willmore, $[\hat{Y}]$ must also be S-Willmore. Each of them could
be obtained as the unique adjoint transform, or the second
envelopping surface of the central sphere congruence, of the other.
$[\hat{Y}]$ is called \emph{the dual Willmore surface} of $[Y]$, and
vice versa. Note that when $\lambda_2\equiv 0$($\lambda_1\equiv 0$),
$[L]$($[R]$) degenerates to a single point. This happens exactly
when $y$ is a null-umbilic surface in $Q^4_1$. Yet the dual Willmore
surface could still be defined if the other $\lambda_i\ne 0$.
\end{remark}
\begin{corollary}
When $y=[Y]:M\rightarrow Q^{4}_{1}$ is a S-Willmore surface without
umbilic points, $[L]$ and $[R]$ are a pair of S-Willmore surfaces
being adjoint transform to each other (one of them might be
degenerate). In particular, the $(-/+)$transforms preserve the
S-Willmore property.
\end{corollary}
\begin{proof}
Since $y$ has no umbilic points, $\lambda_1,\lambda_2$ could not
vanish simultaneously. Without loss of generality, assume
$\lambda_2\ne 0$. Then $[L]$ is an immersion. By the Duality Theorem
above, we see that $[\hat{Y}](=[\tilde{Y}]$) is defined. The
transform chain appeared in Remark~\ref{rem:adjoint} then closes up
as below:
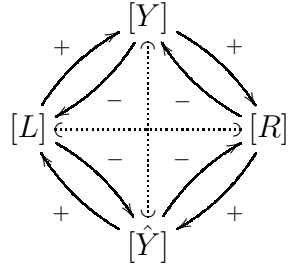
\[
\xymatrix{&&[Y] \ar@/^/[dl]^{-} \ar@/^/[dr]^{+} \ar@{.)}[dd]&\\
&[L] \ar@/^/[ur]^{+} \ar@/^/[dr]^{-} \ar@{.)}[rr]
&&[R]\ar@/^/[ul]^{-} \ar@/^/[dl]^{+} \ar@{.)}[ll]\\
&&[\hat{Y}] \ar@/^/[ul]^{+} \ar@/^/[ur]^{-} \ar@{.)}[uu]& }
\]
It tells us that $[L]$ is the 2-step $(-)$transform and the 2-step
$(+)$transform of $[R]$ at the same time. Equivalently, that means
$[L]$ and $[R]$ are the left and the right adjoint transform of each
other. This proves the conclusion by the Duality Theorem above.
\end{proof}

\section{Spacelike Willmore 2-spheres in $Q^{4}_{1}$}

In this section, we will classify spacelike Willmore 2-spheres in
$Q^{4}_{1}$. This is done by constructing globally defined
holomorphic forms on $S^2$; the vanishing of such forms then enables
us to draw strong conclusions. The reader will see that our method
and result are still similar to the case for Willmore 2-spheres in
$S^4$ [4,16].
\begin{lemma}
(i) Let $y:M\rightarrow Q^{4}_{1}$ be a spacelike Willmore surface
with conformal Hopf differential $\kappa$ for a given coordinate
$z$. Then the 6-form
\begin{equation}
\Theta({\rm d}z)^6=\big[\langle D_{\bar{z}}\kappa,
\kappa\rangle^{2}-\langle \kappa, \kappa\rangle\cdot\langle
D_{\bar{z}}\kappa, D_{\bar{z}}\kappa\rangle\big]({\rm d}z)^6
\end{equation}
is a globally defined holomorphic 6-form on $M$.\\
(ii) When $M=S^2$, we have $\Theta\equiv0$ and $y$ is S-Willmore.
On the subset $M_0\subset M$ where $y$ has no umbilic points, let
$Y$ be the canonical lift of $y$, and $\hat{Y}$ a local lift of
its dual Willmore surface satisfying $\langle
Y,\hat{Y}\rangle=-1$. Then
\begin{equation}
\Omega({\rm d}z)^8=\langle Y_{zz}, Y_{zz}\rangle\langle
\hat{Y}_{zz}, \hat{Y}_{zz}\rangle({\rm d}z)^8
\end{equation}
is a globally defined holomorphic 8-form on $S^{2}.$ So $\Omega
\equiv0$.
\end{lemma}
\begin{proof}
It is easy to verify that these two differential forms are
well-defined (one may use the fact that $\kappa~({\rm
d}z)^\frac{3}{2}({\rm d}\bar{z})^{-\frac12}$ is globally defined).
The holomorphicity of $\Theta({\rm d}z)^{6}$ follows directly from
the Willmore condition \eqref{eq-willmore}.

For conclusion (ii), by the well-known fact that every holomorpic
form on $S^{2}$ must vanish, we know $\Theta\equiv0$. On the other
hand, $\Theta=(\lambda_{1}\gamma_{2}-\lambda_{2}\gamma_{1})^{2}$ by
\eqref{eq-kappa}\eqref{eq-Dkappa}. So on $S^2$ we have
$\lambda_{1}\gamma_{2}-\lambda_{2}\gamma_{1}=0.$ It is just the
S-Willmore condition. Thus on $M_0$ where $\kappa\ne 0$, there is
$D_{\bar{z}}\kappa=-\frac{\bar\mu}{2}\kappa$ for some local function
$\mu$. Define $\hat{Y}$ and $\rho$ as in \eqref{eq-yhat} and
\eqref{eq-rho}, and compute $\hat{Y}_{zz}$ using \eqref{eq-moving2}.
We get $\langle \hat{Y}_{zz},
\hat{Y}_{zz}\rangle=-2\rho^2\lambda_1\lambda_2$. Hence
\[
\Omega=\langle Y_{zz}, Y_{zz}\rangle\langle \hat{Y}_{zz},
\hat{Y}_{zz}\rangle=4(\rho\lambda_1\lambda_2)^2.
\]
Note that in the S-Willmore case
\begin{equation}\label{eq-sigma0}
\sigma=\bar\mu\lambda_1+2\gamma_1=0 =\bar\mu\lambda_2+2\gamma_2.
\end{equation}
So $\rho_{\bar{z}}=\bar\mu\rho$ according to \eqref{eq-Drho}. On the
other hand,
\[
(\lambda_1\lambda_2)_{\bar{z}}
=-\frac{1}{2}\langle\kappa,\kappa\rangle_{\bar{z}} =-\langle
D_{\bar{z}}\kappa,\kappa\rangle = \frac
{\bar\mu}{2}\langle\kappa,\kappa\rangle =-\bar\mu\lambda_1\lambda_2.
\]
Combined together, they show that
$(\rho\lambda_1\lambda_2)_{\bar{z}}=0$ and $\Omega({\rm d}z)^8$ is a
holomorphic differential form defined on $M_0$.

To show $\Omega({\rm d}z)^8$ extends to $M$ as a holomorphic form,
note that by \eqref{eq-sigma0},
\[
\bar\mu_z \lambda_1\lambda_2
=(\bar\mu\lambda_1\lambda_2)_z-\bar\mu(\lambda_1\lambda_2)_z
=(-2\gamma_1\lambda_2)_z
+2\gamma_1(\lambda_2)_z+2\gamma_2(\lambda_1)_z
\]
is a smooth function (depending on $z$). Then for
$\rho=\bar\mu_z+2\lambda_1\bar{\lambda}_2
+2\lambda_2\bar{\lambda}_1$ \eqref{eq-rho}, we see that
$(\rho\lambda_1\lambda_2)^2({\rm d}z)^8$ extends smoothly to $M$
as desired. It is holomorphic both on $M_0$ and in the interior of
$M\setminus M_0$ (it vanishes in the latter case). So it is
holomorphic on the whole $M=S^2$. This completes the proof.
\end{proof}

\begin{theorem}\label{thm-2sphere}
Let $y:S^2\rightarrow Q^4_1$ be a spacelike Willmore 2-sphere. Then
it must be a surface among the following two classes:

(i) it is the conformal compactification of a stationary surface in
$\mathbb{R}^4_1$.

(ii) it is one of the polar surfaces of a surface in class (i).
\end{theorem}
\begin{proof}
First we observe that \eqref{eq-gamma}\eqref{eq-willmore3} may be
re-written as
\[
\left\{\begin {array}{ll}
{\lambda_1}_{\bar{z}}=-\bar\alpha\lambda_1+\gamma_1,\\
{\gamma_1}_{\bar{z}}=-\frac{\bar{s}}{2}\lambda_1-\bar\alpha\gamma_1.
\end {array}\right.
\]
By a lemma of Chern (see section 4 in [7]), either $\lambda_1$ is
identically zero on $S^2$, or it has only isolated zeroes. The same
conclusion holds for $\lambda_2$. Now that we have shown
$\rho\lambda_{1}\lambda_{2}\equiv 0$, one of $\rho,\lambda_1,
\lambda_2$ must vanish identically on $S^2$.

If $\rho\equiv0$, $[\hat{Y}]$ degenerates to a single point due to
\eqref{eq-Dyhat} and $\sigma=0$. Applying a transformation $T\in
O(4,2)$ if necessary, we can set $\hat{Y}=(1,0,0,0,0,1)$ and
$Y=(\frac{-1+\langle u,u\rangle}{2},u,\frac{1+\langle
u,u\rangle}{2})$ where $u:U\rightarrow\mathbb{R}^{4}_{1}$. Let $z$
be an arbitrary complex coordinate. Then we have
\[
Y_{z\bar{z}}=a Y+\langle Y_{z},Y_{\bar{z}}\rangle N,
~~\hat{Y}=N+\bar{\mu}Y_{z}+\mu Y_{\bar{z}}+\langle
Y_{z},Y_{\bar{z}}\rangle|\mu|^{2}Y ,
\]
where $a,\mu$ are two functions. It is easy to see$$
\hat{Y}_{z}=-\mu\langle Y_{z},Y_{\bar{z}}\rangle\hat{Y}+\cdots.$$ So
$\mu\equiv0$ and $Y_{z\bar{z}}=aY+\langle
Y_{z},Y_{\bar{z}}\rangle\hat{Y}$. Replacing by $u$ leads to
$$(\langle u_{z\bar{z}},u\rangle,u_{z\bar{z}},\langle u_{z\bar{z}},u\rangle)=(\frac{-a+a\langle u,u\rangle}{2},au,\frac{a+a\langle u,u\rangle}{2}).$$
This implies $a\equiv0$ and $u_{z\bar{z}}\equiv0$. So $u$ is a
stationary surface in $\mathbb{R}^4_1$, and $y=[Y]$ belongs to
class (i). \footnote{An alternative proof is by the meaning of the
mean curvature sphere. Since every central sphere of $y$ passing
through a fixed point $[\hat{Y}]$ of $Q^4_1$, which could be
viewed as a point at infinity for some affine $\mathbb{R}^4_1$,
each sphere is a plane in this $\mathbb{R}^4_1$.
Corollary~\ref{cor:mini} implies the conclusion.}

If $\lambda_1\equiv0$($\lambda_2\equiv0$), $[R]$($[L]$) is a point.
Using the conclusion in (i) for surface $[L]$($[R]$) finishes the
proof.
\end{proof}
\begin{remark}
Note that surfaces of class (ii) are exactly spacelike null-umbilic
Willmore surfaces. So one has
\[
\langle\kappa,\bar\kappa\rangle=
-\lambda_1\bar{\lambda}_2-\lambda_2\bar{\lambda}_1=0.
\]
As a consequence, its induced conformal metric
$\langle\kappa,\bar\kappa\rangle({\rm d}z)^2$ as well as the
Willmore functional is always zero, which is different from the
case in $S^{4}$. Here a left question is: \emph{For a space-like
Willmore 2-sphere in $Q^{4}_{1}$, if its Willmore functional
equals zero, must it be of type (ii)?}
\end{remark}

\section{Examples}

First let us see some special Willmore surfaces contained in a
3-dimensional space.
\begin{example}
Embed $\mathbb{R}^{3}\subset\mathbb{R}^{4}_{1}$ via $u\to (u,1)$,
$\mathbb{R}^{3}_{1}\subset\mathbb{R}^{4}_{1}$ via $u\to (1,u)$.

(i) Let $u:M^{2}\rightarrow\mathbb{R}^{3}$ be a minimal surface.
Then $(u,1):M^{2}\rightarrow\mathbb{R}^{4}_{1}$ is a spacelike
stationary surface in $\mathbb{R}^{4}_{1}$, and
\[
Y=\Big(\frac{\langle u,u\rangle}{2}-1,u,1,\frac{\langle
u,u\rangle}{2}\Big) :M^2 \rightarrow C^5
\]
gives a spacelike S-Willmore surface $[Y]:M^{2}\rightarrow
Q^{4}_{1}$. (Essentially this comes from the conformal embedding
$\varphi_0$.) Let $g:M^{2}\rightarrow S^{4}_{1}\subset R^{5}_{1}$
denote the conformal Gauss map of $u$ as in [2] and
$e=(-1,0,0,0,-1,1)\in\mathbb{R}^{6}_{2}$. It is straightforward to
verify that $[e+(g,0)]$ and $[e-(g,0)]$ are the polar surfaces of
$[Y]$.

(ii) Let $u:M^{2}\rightarrow\mathbb{R}^{3}_{1}$ be a spacelike
maximal surface. Then $(1,u):M^{2}\rightarrow\mathbb{R}^{4}_{1}$ is
a spacelike stationary surface in $\mathbb{R}^{4}_{1}$, and
\[
Y=\Big(\frac{\langle u,u\rangle}{2},1,u,\frac{\langle
u,u\rangle}{2}+1\Big):M^{2}\rightarrow C^5
\]
gives a spacelike S-Willmore surface $[Y]:M^{2}\rightarrow
Q^{4}_{1}$. Let $\tilde{g}:M^{2}\rightarrow H^{4}_{1}\subset
R^{5}_{2}$ denote the conformal Gauss map of $u$ as in $[1]$ and
$\tilde{e}=(1,1,0,0,0,1)\in\mathbb{R}^{6}_{2}$. Then
$[\tilde{e}+(0,\tilde{g})]$ and $[\tilde{e}-(0,\tilde{g})]$ are the
polar surfaces of $[Y]$.

(iii) Suppose $u:M^{2}\rightarrow\mathbb{R}^{3}$ is a Laguerre
minimal surface and $n:M^{2}\rightarrow S^{2}$ its Gauss map. Its
Laguerre lift
\[
Y=(n,u\cdot n,-u\cdot n,1):M^{2}\rightarrow C^5
\]
gives a spacelike S-Willmore surface $[Y]:M^{2}\rightarrow
Q^{4}_{1}$. We denote $g': M^{2}\rightarrow R^{4}_{1}\hookrightarrow
C^5$ its Laguerre Gauss map (see [18]). Then the point
$[(0,0,0,1,-1,0)]$ and $[g']$ are the polar surfaces of $[Y]$.
\end{example}

\begin{example}
Consider a spacelike Willmore surface $y=[Y]$ both of type $(i)$ and
type $(ii)$ as in Theorem~\ref{thm-2sphere}. That means either of
$[\hat{Y}]$ and $[L]$ is a single point, and $[Y]$ is the conformal
compactification of a stationary surface
$x:M^{2}\rightarrow\mathbb{R}^{4}_{1}$. Without loss of generality,
suppose
\[
\hat{Y}=(1,0,0,0,0,1),~L=(0,1,0,0,1,0),~~ Y=\Big(\frac{-1+\langle
x,x\rangle}{2},x,\frac{1+\langle x,x\rangle}{2}\Big).
\]
From $\langle Y,L\rangle=0$, we see that the surface
$x=(x_{1},x_{2},x_{3},x_{4})$ must satisfy $x_{1}=x_{4}$, which
means that $x$ in fact is a zero mean curvature surface in
$\mathbb{R}^{3}_{0}\subset\mathbb{R}^{4}_{1}$. For details of such
surfaces, see [18].
\end{example}

Among compact surfaces, 2-spheres and tori are simplest and most
important. In general, Willmore tori are not necessarily S-Willmore
surfaces. Here we give such a class of spacelike Willmore tori which
are homogenous.

\begin{example}
Let$$e_{1}=(\cos\frac{t\theta}{\sqrt{t^{2}-1}}\cos\phi,
\cos\frac{t\theta}{\sqrt{t^{2}-1}}\sin\phi,
\sin\frac{t\theta}{\sqrt{t^{2}-1}}\cos\phi,
\sin\frac{t\theta}{\sqrt{t^{2}-1}}\sin\phi),$$
$$e_{2}=\frac{\partial e_{1}}{\partial\phi}=e_{1\phi},
e_{3}=\frac{\sqrt{t^{2}-1}}{t}e_{1\theta},
e_{4}=\frac{\sqrt{t^{2}-1}}{t}e_{2\theta},$$
where $t>1$.
 Let
\begin{equation} \label{eq:1}
\begin{split}
Y_{t}(\theta,\phi)&: \mathbb{R}\times\mathbb{R}
\longrightarrow \mathbb{R}^{6}_{2} \\
 Y_{t}(\theta,\phi)&=(e_{1},\cos\frac{\theta}{\sqrt{t^{2}-1}},
 \sin\frac{\theta}{\sqrt{t^{2}-1}}).
 \end{split}
 \end{equation}
For simplicity, we omit the subscript $"_{t}"$ of $Y_{t}$. We have
that $y=[Y]:\mathbb{R}\times\mathbb{R}\rightarrow Q^{4}_{1}$ is a
spacelike Willmore torus of $Q^{4}_{1}$ for any rational number
$t>1$.

For the lift $Y$ we set $z=\theta+i\phi$. It is easy to verify $Y$
is a canonical lift with respect to $z$. We have
\begin{equation} \label{eq:5}
\left\{ \begin{aligned}
         Y_{z\bar{z}} &= -\frac{t^{2}}{4(t^{2}-1)}Y+\frac{1}{2}N, \\
         Y_{zz} &=  -\frac{1}{4(t^{2}-1)}Y-\frac{it}{2\sqrt{2}\sqrt{t^{2}-1}}(L-R),\\
         L_{z} &=  -\frac{i}{2\sqrt{t^{2}-1}}L-\frac{t}{2\sqrt{2}(t^{2}-1)}Y+\frac{it}{\sqrt{2}\sqrt{t^{2}-1}}Y_{\bar{z}}, \\
         R_{z} &= \frac{i}{2\sqrt{t^{2}-1}}R-\frac{t}{2\sqrt{2}(t^{2}-1)}Y-\frac{it}{\sqrt{2}\sqrt{t^{2}-1}}Y_{\bar{z}},  \\
                  N_{z}&= -\frac{t^{2}}{2(t^{2}-1)}Y_{z}-\frac{1}{2(t^{2}-1)}Y_{\bar{z}}+\frac{t}{2\sqrt{2}(t^{2}-1)}L+\frac{t}{2\sqrt{2}(t^{2}-1)}R.                         \end{aligned} \right.
                          \end{equation}

Here
\begin{equation} \label{eq:4}
\left\{ \begin{aligned}
         Y_{z} &=\frac{1}{2\sqrt{t^{2}-1}}(te_{3}-i\sqrt{t^{2}-1}e_{2},-\sin\frac{\theta}{\sqrt{t^{2}-1}},\cos\frac{\theta}{\sqrt{t^{2}-1}}),\\
         N &= \frac{1}{2}(-e_{1},\cos\frac{\theta}{\sqrt{t^{2}-1}},\sin\frac{\theta}{\sqrt{t^{2}-1}}), \\
         L &= \frac{1}{\sqrt{2}\sqrt{t^{2}-1}}(\sqrt{t^{2}-1}e_{4}+e_{3},-t\sin\frac{\theta}{\sqrt{t^{2}-1}},t\cos\frac{\theta}{\sqrt{t^{2}-1}}),  \\
         R &= \frac{1}{\sqrt{2}\sqrt{t^{2}-1}}(-\sqrt{t^{2}-1}e_{4}+e_{3},-t\sin\frac{\theta}{\sqrt{t^{2}-1}},t\cos\frac{\theta}{\sqrt{t^{2}-1}}).                         \end{aligned} \right.
                          \end{equation}

So it is easy to see that $Y$ is spacelike Willmore and not
S-Willmore.

The adjoint surface of $Y$ with respect to $L$ is
$$\hat{Y}=\frac{1}{2}\left(\frac{2-t^{2}}{t^{2}-1}e_{1}+\frac{1}{\sqrt{t^{2}-1}}e_{2},\frac{t^{2}}{t^{2}-1}\cos\frac{\theta}{\sqrt{t^{2}-1}},\frac{t^{2}}{t^{2}-1}\sin\frac{\theta}{\sqrt{t^{2}-1}}\right).$$

The adjoint surface of $Y$ with respect to $R$ is
$$\tilde{Y}=\frac{1}{2}\left(\frac{2-t^{2}}{t^{2}-1}e_{1}-\frac{1}{\sqrt{t^{2}-1}}e_{2},\frac{t^{2}}{t^{2}-1}\cos\frac{\theta}{\sqrt{t^{2}-1}},\frac{t^{2}}{t^{2}-1}\sin\frac{\theta}{\sqrt{t^{2}-1}}\right).$$

We point out that $Y$ is a homogenous torus which is the orbit of
the subgroup
$$G=\left(
      \begin{array}{cccccc}
        e_{1}^{T} & e_{2}^{T} & e_{3}^{T} & e_{4}^{T} & 0 & 0 \\
        0 & 0 & 0 & 0 &\cos\frac{\theta}{\sqrt{t^{2}-1}} &  -\sin\frac{\theta}{\sqrt{t^{2}-1}} \\
        0 & 0 & 0 & 0 &\sin\frac{\theta}{\sqrt{t^{2}-1}} & \cos\frac{\theta}{\sqrt{t^{2}-1}}
      \end{array}
    \right)
$$
acting on $(1,0,0,0,1,0)^{T}$. Here $^{T}$ denotes transposition.

If $y_{t}$ is a torus, then $t$ must be some rational number.
Suppose $t=\frac{p}{q}$, where $p,q\in \mathbb{N}$. Then the
Willmore functional of $y_{t}$ is
\begin{equation}
W(y_{t})=\frac{p^{2}}{\sqrt{p^{2}-q^{2}}}\pi^{2}
\end{equation}
So the minimum of Willmore functional of $y_{t}$ is
$\frac{4}{\sqrt{3}}\pi^{2}$.
\end{example}

\def\refname{Reference}

\end{document}